\documentclass[preprint, 12pt, draft]{elsarticle}

\usepackage{amsthm,amsmath,amsfonts,amssymb,latexsym, mathrsfs}

\theoremstyle{plain}
\newtheorem{theorem}{Theorem}

\numberwithin{corollary}{theorem}

\journal{Statistics \& Probability Letters}

\begin{document}

\setcitestyle{authoryear, round, comma, aysep={;}, yysep={,}, notesep={,}}

\begin{frontmatter}

\title{Integral form of the COM--Poisson normalization constant}

\author{Tibor K. Pog\'any}

\address{Faculty of Maritime Studies, University of Rijeka, 51000 Rijeka, Studentska 2, Croatia {\it and} \\ 
Institute of Applied Mathematics, \'Obuda University, 1034 Budapest, B\'ecsi \'ut 96/b, Hungary} 

\ead{poganj@pfri.hr} 
\ead[url]{http://www.pfri.hr/~poganj}

\begin{abstract}
In this brief note an integral expression is presented for the COM--Poisson renormalization constant $Z(\lambda, \nu)$ 
on the real axis. 
\end{abstract}

\end{frontmatter}

\allowdisplaybreaks

\section{Introduction and Motivation} 

Consider a rv $X$ having classical COM--Poisson distribution initiated by \citet{Conway} in the following form:
   \[ \mathsf P(X=n) = \frac{\lambda^n}{Z(\lambda, \nu)\, (n!)^\nu}, \qquad \lambda, \nu>0; \quad n\in \mathbb N_0,\]
$Z(\lambda, \nu)$ being the reciprocal of the normalizing constant. 

To the best of my knowledge, the only existing integral expression for $Z(\lambda, \nu)$ was given by \citet[p. 141, Eq. (40)]{SMKBB} 
for positive integer $\nu$ as a multiple integral 
   \begin{align*}
	    Z(\lambda, \nu) &= \frac1{(2\pi)^{\nu-1}} \int_{-\pi}^\pi \cdots \int_{-\pi}^\pi \exp\Bigg\{ \sum_{j=1}^{\nu-1} 
	                       \exp \big\{{\rm i}x_j\big\} \\ 
											&\qquad \qquad + \lambda \exp\Big\{-{\rm i} \sum_{j=1}^{\nu-1}x_j\Big\} \Bigg\}\,
												 {\rm d}x_1 \cdots {\rm d}x_{\nu-1}\, .
	 \end{align*}
In the same case with $\nu \geq 1$ assumed to be an integer, Nadarajah reported on the closed form expression \citet[p. 619, Eq. (6)]{Nad}  
   \[ Z(\lambda, \nu) = {}_0F_\nu(-; 1, \cdots, 1; \lambda)\,,\]
where ${}_0F_\nu$ stands for the familiar generalized hypergeometric function with zero upper and $\nu$ lower parameters. 

Another fashion related approximation and extension results regarding $Z(\lambda, \nu)$ are numerous, see for instance 
the recent articles \citet{SMKBB, Gillispie, Simsek} and the references therein. 

\section{The integral expression for $Z(\lambda, \nu)$}

Here a single definite integral expression is established for $Z(\lambda, \nu)$, for all positive $\nu$. 

\begin{theorem} For all $\lambda, \nu>0$ we have
   \[ Z(\lambda, \nu) = \frac1{1-\lambda} + \frac\nu{\lambda-1} \int_1^\infty x^{-(\nu+1)} \, 
											   \lambda^{\left[ \Gamma^{-1}(x)\right]}\, {\rm d}x\,,\]
{\it where $\Gamma^{-1}$ stands for the increasing branch of the inverse Gamma in the right half-plane, while $[x]$ denotes the 
integer part of some real $x$.} \bigskip
\end{theorem} 

\noindent {\sc Proof.}
It turns out that
   \[ Z(\lambda, \nu) = \sum_{n\geq 0} \frac{\lambda^n}{(n!)^\nu} = \sum_{n\geq 0} \lambda^n\, {\rm e}^{-\nu \ln \Gamma(n+1)}\, ,\]
is a classical Dirichlet series. The main tool we refer to is the Cahen formula for the Laplace 
integral representation of Dirichlet series (reported firstly without proof by \citet{Cah} and proved by \citet{Perron}). 
Namely, the Dirichlet series
   \[ \mathscr D_{\boldsymbol a}(r) = \sum_{n \ge 1} a_n e^{-r b_n},\]
where $\Re(r)>0,$ possessing positive monotone increasing divergent to infinity sequence $(b_n)_{n \ge 1}$, has a Laplace   
integral representation \citet[p. 97]{Cah}
   \[ \mathscr D_{\boldsymbol a}(r) = r \int_0^\infty e^{-rx} \sum_{n \colon b_n \leq x} a_n\, {\rm d}x 
			                              = r \int_0^\infty e^{-rx} \sum_{n=1}^{[b^{-1}(x)]} a_n\, {\rm d}x \, ,\]
since $b \colon \mathbb R_+ \mapsto \mathbb R_+$ is monotone, and there exists unique inverse $b^{-1}$ for the function 
$b \colon \mathbb R_+ \mapsto \mathbb R_+$, being $b|_{\mathbb N} = (b_n)$, see also \citet{P2}. 

Thus, bearing in mind that $b(x) \equiv \ln \Gamma(x+2)$ is monotone and invertible for $x \geq 1$, being $x> \alpha = 1.4616 \cdots$  
the abscissa of the minimum of the Gamma function \citet{Ped}, we only have to follow the previous derivation procedure:
   \begin{align*}
	    Z(\lambda, \nu) &= \sum_{n\geq 0} \lambda^n\, {\rm e}^{-\nu \ln \Gamma(n+1)} \\
			                &= 1+\lambda + \lambda^2 \nu \int_0^\infty {\rm e}^{-\nu x} \, 
											   \sum_{n \colon \ln \Gamma(n+2) \leq x} \lambda^{n-1}\, {\rm d}x \\
											&= 1+\lambda + \lambda^2 \nu \int_0^\infty {\rm e}^{-\nu x} \, 
											   \sum_{n = 1}^{\left[ \Gamma^{-1}({\rm e}^x)\right]-2} \lambda^{n-1}\, {\rm d}x \\
										  &= 1+\lambda + \frac{\lambda^2 \nu}{\lambda-1} \int_0^\infty {\rm e}^{-\nu x} \, 
											   \left(\lambda^{\left[ \Gamma^{-1}({\rm e}^x)\right]-2} -1\right) \, {\rm d}x \\
											&= \frac1{1-\lambda} + \frac\nu{\lambda-1} \int_1^\infty x^{-(\nu+1)} \, 
											   \lambda^{\left[ \Gamma^{-1}(x)\right]}\, {\rm d}x\,, 
	 \end{align*} 
which finishes the proof of the stated integral representation.  \hfill $\Box$

\section*{Acknowledgement} 
I am grateful to Professor Satish Iyengar, University of Pittsburgh for discussion in which he introduced me the 
COM--Poisson distribution's open questions.


\begin{thebibliography}{7}

\bibitem[Cahen(1894)]{Cah} 
Cahen, E. (1894). Sur la fonction $\zeta(s)$ de Riemann et sur des fontions analogues. 
{\it Ann. Sci. l'\'Ecole Norm. Sup. S\'er. Math.} {\bf 11}, 75--164.

\bibitem[Conway--Maxwell(1962)]{Conway} 
Conway, R. W. and Maxwell, W. L. (1962). A queuing model with state dependent service rates. {\it J. Industrial Engineering} 
{\bf 12}, 132--136.

\bibitem[Gillispie--Christopher(2015)]{Gillispie}
Gillispie, S. B. and Christopher, G. C. (2015). Approximating the Conway-Maxwell-Poisson distribution normalization constant. 
{\it Statistics} {\bf 49}, No. 5, 1062--1073. 

\bibitem[Nadarajah(2009)]{Nad} 
Nadarajah, S. (2009). Useful moment and CDF formulations for the COM–Poisson distribution. {\it Stat. Papers} {\bf 50}, 617--622.

\bibitem[Pedersen(2015)]{Ped} 
Pedersen, H. L. (2015). Inverses of gamma functions. {\it Constr. Approx.} {\bf 41}, No. 2, 251--267.

\bibitem[Perron(1908)]{Perron} 
Perron, O. X. (1908). Zur Theorie der Dirichletschen Reihen. {\em J. Reine Angew. Math.} {\bf 134}, 95--143. 

\bibitem[Pog\'any(2005)]{P2} 
Pog\'any, T. K. (2005). Integral representation of Mathieu $(\boldsymbol a, \boldsymbol \lambda)$--series.  
{\it Integral Transforms Spec. Funct.} {\bf 16}, No.8, 685--689. 

\bibitem[Schmueli {\it et al.}(2005)]{SMKBB} 
Shmueli, G., Minka, T. P.,  Kadane,J. B., Borle, S. and Boatwright, P. (2005). A useful distribution for fitting discrete
data: revival of the Conway-Maxwell-Poisson distribution. {\it J. Roy. Statist. Soc. Ser. C} {\bf 54} No. 1, 127--142. 

\bibitem[\c Sim\c sek--Iyengar(2015)]{Simsek} 
\c Sim\c sek, B. and S. Iyengar, S. (2015). Approximating the Conway--Maxwell--Poisson normalizing constant. {\it Filomat} (to appear). 


\end{thebibliography}
\end{document}